\documentclass[12pt,a4paper]{article}

\usepackage[margin=1in]{geometry}
\usepackage{setspace}
\usepackage{times}
\usepackage{graphicx}
\graphicspath{{figures/}}
\usepackage{booktabs}
\usepackage{multirow}
\usepackage{amsmath,amssymb}
\usepackage{array}
\usepackage{caption}
\usepackage{subcaption}
\usepackage{url}
\usepackage{authblk}
\usepackage{hyperref}
\usepackage[numbers,sort&compress]{natbib}
\usepackage{algorithm}
\usepackage{algorithmic}

\title{\textsf{CrackMorph-XAI-Net: A Topology-Preserving and Explainable Framework for Automated Crack Morphology}}
\author[1]{Sri Surya Paravallika Ajjarapu \thanks{sajjarapu@islander.tamucc.edu}}
\author[1,$\dagger$]{S. M. Mallikarjunaiah \thanks{m.muddamallappa@tamucc.edu}}
\affil[1]{Department of Mathematics \& Statistics, Texas A\&M University-Corpus Christi, Texas- 78412, USA}
\affil[$\dagger$]{Corresponding author}
\date{}

\begin{document}

\maketitle

\begin{abstract}
Automated crack inspection is increasingly recognized as a critical component of infrastructure monitoring; however, cracks continue to be reported primarily as binary segmentation masks by many current vision-based systems. While localization is facilitated by such masks, limited structural information is provided for robust engineering interpretation. For practical crack assessment, measurable morphological features—including centerline geometry, branching behavior, junction locations, topology, and severity-related indicators—are required. In this work, \textit{CrackMorph-XAI-Net}, an explainable morphology-aware framework for image-based crack analysis, is presented. Crack image and region-mask data are converted into a sequence of interpretable structural outputs through four distinct stages: topology-preserving skeleton extraction, junction detection via Gaussian heatmap regression, morphology descriptor computation, and severity-oriented screening. To support rigorous stage-wise evaluation, the standard \textit{CRACK500} benchmark is extended with aligned skeleton maps, junction heatmaps, and topology labels. Experimental validation demonstrates that a mean Dice coefficient of 0.991 is achieved by the learned skeleton extraction stage, with topology preserved in 98.5\% of test images. Furthermore, a recall of 0.964 and an F1-score of 0.887 are obtained in the junction detection stage, highlighting the efficacy of heatmap regression for sparse structural targets. Strong agreement between predicted and reference morphology values is revealed by descriptor-level evaluation, with correlations exceeding 0.95 for length, width, orientation, junction count, and tortuosity. When compared against a classical thinning-based structural baseline, significant improvements in skeleton quality, junction recall, descriptor accuracy, and topology classification are exhibited by the proposed framework. Ultimately, these results demonstrate that crack analysis can be effectively advanced beyond mere segmentation through the generation of transparent, measurable, and engineering-relevant structural representations.
\end{abstract}


\section{Introduction}
\label{sec:introduction}

Cracks are among the most visible early indicators of deterioration in pavements, bridges, concrete surfaces, and broader infrastructure systems. Arising from repeated loading, aging, environmental exposure, material defects, or construction-related factors, these anomalies demand early detection and characterization. Left unaddressed, initial surface cracks can propagate into severe structural damage, escalating maintenance costs, diminishing service life, and introducing significant safety risks. Consequently, rigorous crack inspection plays an indispensable role in structural health monitoring and preventative maintenance planning.

Traditionally, crack inspection has been performed manually or augmented with semi-automated tools. While manual inspection benefits from direct human judgment, it is inherently time-consuming, labor-intensive, and subjective. Discrepancies in severity ratings frequently occur among inspectors due to variations in experience, ambient lighting conditions, surface texture, and contextual field variables. Across large-scale transportation and infrastructure networks, these limitations preclude consistent and scalable assessments. To address these inefficiencies and mitigate subjective variation, automated computer vision methodologies have been continuously developed over the past few decades.

Early automated crack detection methods relied heavily on classical image processing operations, such as thresholding, morphological filtering, and edge detection using Sobel, Canny, or Laplacian filters to highlight intensity discontinuities \citep{abdel2003edge}. While computationally efficient, these operators proved highly sensitive to shadows, stains, non-uniform illumination, and complex aggregate patterns, necessitating extensive post-processing. Subsequent graph-based and minimum-path approaches, such as CrackTree \citep{zou2012cracktree}, sought to improve continuity by representing cracks as connected structures rather than isolated dark pixels. Concurrently, feature-based machine learning methods extracted handcrafted descriptors—including local binary patterns, Gabor features, and histogram of oriented gradients—to train classifiers like support vector machines and random forests \citep{cord2012road,oliveira2013road,shi2016random}. Although these techniques offered greater adaptability than fixed thresholding, their success remained rigidly tied to optimal descriptor design, and their primary output was constrained to basic detection or region segmentation rather than detailed structural morphology.

The landscape of automated crack analysis was fundamentally transformed by deep learning, which enabled robust hierarchical feature extraction directly from raw imagery. Early convolutional neural networks for crack damage detection \citep{cha2017deep} quickly evolved into fully convolutional networks capable of dense pixel-level prediction \citep{long2015fcn}. Architectures like U-Net \citep{ronneberger2015unet} and DeepCrack \citep{liu2019deepcrack} became highly effective at preserving fine spatial details through skip connections and multi-scale feature fusion, addressing the challenge of segmenting thin, elongated crack structures. More recently, transformer-based models like CrackFormer \citep{liu2022crackformer} and advanced attention mechanisms \citep{yang2025pcmsda} have introduced long-range dependency modeling to further refine crack localization under challenging visual conditions.

Despite these significant advancements in pixel-level localization, the dominant output of most contemporary systems remains a binary or probabilistic segmentation mask. While a mask effectively localizes damage, it fails to capture how a crack is structurally organized. For rigorous engineering interpretation, a mask alone is insufficient; maintenance decisions are dictated not just by the presence of a defect, but by quantitative morphological descriptors such as length, width, orientation, branching behavior, tortuosity, connectivity, and topological classification. Historically, these metrics have been approximated post-segmentation using classical skeletonization or medial-axis transforms. However, such heuristic post-processing is notoriously brittle, with minor boundary noise in a mask frequently generating false branches, broken centerlines, or distorted topology \citep{mosinska2018topology,zhang2018road}. This highlights a critical need to transition from classical post-processing to learning-based structural prediction.

Furthermore, as automated systems are integrated into safety-critical engineering workflows, explainable artificial intelligence (XAI) becomes essential. End-users require transparency to trust algorithmic assessments \citep{arrieta2020xai}. In the realm of fracture mechanics and fatigue crack analysis, post-hoc explanation methods like Grad-CAM have been utilized to verify that models attend to physically meaningful phenomena, such as crack-tip fields, rather than image artifacts \citep{melching2022xai,roux2009dic,zhao2019dic,rethore2015cracktip}. However, post-hoc visualizations layered over black-box segmentation models do not inherently yield measurable, structural data. 

To bridge the gap between pixel-level localization and transparent, morphology-aware interpretation, this paper introduces CrackMorph-XAI-Net. Operating under the assumption that a crack image and its corresponding region mask are available, the proposed framework explicitly converts raw imagery into an interpretable sequence of structural representations. Rather than appending explainability as an afterthought, the framework is explainable by design: it sequentially predicts crack skeletons, detects junctions, computes morphology descriptors, and derives topology classes. Each intermediate product serves as a measurable, verifiable foundation for the final severity-oriented screening output. To address these limitations, this work proposes a comprehensive morphology-aware crack analysis framework designed to advance automated inspection beyond conventional binary segmentation and toward rigorous, engineering-relevant structural interpretation. Central to this framework is the introduction of a learned skeleton extraction stage, which recovers topology-preserving crack centerlines and significantly curtails the noise, broken connectivity, and structural artifacts typically associated with classical thinning algorithms. Building upon this refined skeletal foundation, a Gaussian heatmap regression strategy is implemented to achieve precise junction detection, effectively overcoming the inherent optimization challenges associated with learning highly sparse branch-point targets. These predictive stages culminate in a deterministic descriptor formulation that computes essential engineering metrics directly from the generated structural outputs, yielding transparent measurements for crack length, width, orientation, junction count, tortuosity, and overall topology class. Finally, to facilitate rigorous, stage-wise evaluation and support future benchmarking of these specialized techniques, this study develops a morphology-aware extension of the widely used CRACK500 dataset, newly featuring fully aligned skeleton maps, junction heatmaps, and comprehensive topology labels.

The remainder of this paper is organized as follows. Section~\ref{sec:dataset} describes the dataset extension and preparation. Section~\ref{sec:method} details the architecture and operational stages of the proposed framework. Section~\ref{sec:experiments} outlines the experimental setup, training protocols, and evaluation metrics. Section~\ref{sec:results} presents the quantitative and qualitative results of the study. Section~\ref{sec:discussion} explores the interpretability, practical implications, and limitations of the approach, and Section~\ref{sec:conclusion} concludes the paper.

\section{Dataset}
\label{sec:dataset}

\subsection{Base Benchmark}

This study leverages the widely recognized CRACK500 dataset \citep{yang2020crack500} as its foundational benchmark. Originally curated to advance pavement crack segmentation, CRACK500 comprises 500 high-resolution images paired with meticulously annotated pixel-level binary crack masks. To ensure rigorous comparative evaluation, the standard dataset split is maintained: 250 images are allocated for training, 50 for validation, and 200 for independent testing. The dataset is particularly advantageous due to its inclusion of diverse pavement textures, complex crack geometries, fluctuating illumination conditions, and challenging background noise, making it a highly robust benchmark for automated segmentation research.

However, while CRACK500 excels in evaluating binary localization architectures, it was not inherently designed to facilitate morphology-aware structural learning. The provided binary masks successfully isolate damaged regions but completely lack the explicit encoding required for advanced engineering interpretation—namely, underlying centerlines, branching behavior, junction coordinates, and macroscopic topology classes. Consequently, to accurately evaluate crack morphology and structural geometry beyond mere pixel classification, augmented structural annotations are strictly necessitated.

\subsection{CRACK500-Extended}

To explicitly address these limitations and provide the necessary supervisory signals for the proposed framework, the base CRACK500 dataset is comprehensively expanded into a morphology-aware benchmark, denoted herein as CRACK500-Extended. While retaining the original input images and corresponding binary masks, this extended dataset introduces three critical, spatially aligned structural annotation layers:
\begin{enumerate}
    \item \textbf{Skeleton maps:} One-pixel-wide, topology-preserving centerline representations of the crack regions.
    \item \textbf{Junction heatmaps:} Continuous, two-dimensional Gaussian distributions centered precisely at annotated crack branch points.
    \item \textbf{Topology labels:} Global, image-level complexity categories derived from the localized branching structure.
\end{enumerate}

The generated skeleton maps serve as the direct supervisory target for the learned centerline recovery stage. Initial skeleton candidates are algorithmically derived from the base binary masks via medial-axis transformations, followed by rigorous refinement protocols to eliminate spurious micro-branches, boundary-induced artifacts, and artefactual discontinuities. Subsequently, junction annotations are deterministically extracted from these refined skeletal structures. Crucially, rather than representing each structural branch point as an isolated, single-pixel target—which penalizes minor spatial deviations and exacerbates class imbalance—junction locations are transformed into Gaussian heatmaps. This continuous probabilistic representation provides a smoother optimization landscape and richer spatial gradients during network training. Finally, the topology labels abstract the underlying structural data into highly interpretable global categories, classifying the overall crack geometry as either linear, branched, complex, or a full network.

\begin{figure}[t]
\centering
\includegraphics[width=0.98\linewidth]{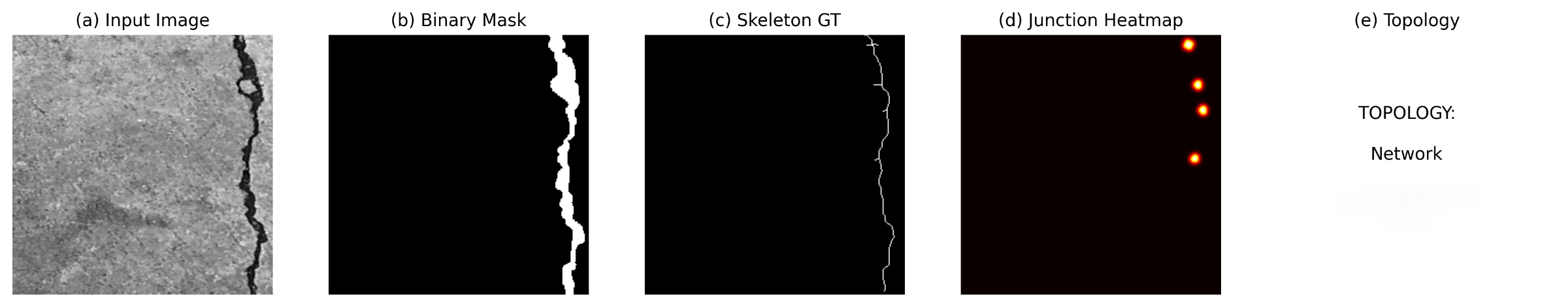}
\caption{CRACK500-Extended annotation layers utilized for morphology-aware crack analysis, sequentially illustrating the input image, binary crack mask, skeleton representation, junction heatmap, and topology label.}
\label{fig:dataset_layers}
\end{figure}

\begin{table}[t]
\centering
\caption{Summary comparison of the base and extended CRACK500 annotation layers.}
\label{tab:dataset}
\begin{tabular}{lcccccc}
\toprule
Dataset & Images & Split & Mask & Skeleton & Junction heatmap & Topology label \\
\midrule
CRACK500 & 500 & 250/50/200 & Yes & No & No & No \\
CRACK500-Extended & 500 & 250/50/200 & Yes & Yes & Yes & Yes \\
\bottomrule
\end{tabular}
\end{table}

\subsection{Preprocessing and Augmentation}

Prior to network ingestion, a standardized preprocessing pipeline is applied to ensure dimensional consistency and optimal model convergence. All input images are uniformly resized to a resolution of $640 \times 640$ pixels, converted to a single-channel grayscale format to minimize unnecessary computational overhead, and normalized to a floating-point range of $[0, 1]$. To strictly preserve label discreteness and structural integrity, the corresponding binary masks and discrete skeleton maps are resized utilizing nearest-neighbor interpolation. Conversely, bilinear interpolation is deliberately applied to the junction heatmaps to maintain their continuous Gaussian gradient structure.

To mitigate overfitting and enhance the model's geometric invariance, the training set is subjected to a rigorous spatial augmentation regimen. This includes horizontal and vertical flipping, arbitrary rotations, and elastic deformations designed to simulate complex perspective distortions encountered in field inspections. Crucially, all spatial transformations are applied synchronously across the raw image, mask, skeleton, and heatmap tensors to guarantee perfect structural alignment. Following the exclusion of augmented samples exhibiting severe boundary cropping artifacts, the final augmented training repository is expanded to 1,896 highly diverse instances.

\section{Proposed Method}
\label{sec:method}

\subsection{Framework Overview}

The proposed framework is conceptualized as a transparent, stage-wise pipeline explicitly tailored for morphology-aware crack analysis. Departing from conventional automated paradigms that yield structurally opaque binary segmentation masks, this architecture systematically transforms pixel-level crack region data into quantifiable structural representations and robust geometric descriptors. The overarching pipeline is modularized into four sequential, fully interpretable stages: 
\begin{enumerate}
    \item topology-preserving skeleton extraction,
    \item junction detection utilizing Gaussian heatmap regression,
    \item deterministic morphology descriptor computation, and
    \item severity-oriented screening.
\end{enumerate}

\begin{figure}[t]
\centering
\includegraphics[width=0.98\linewidth]{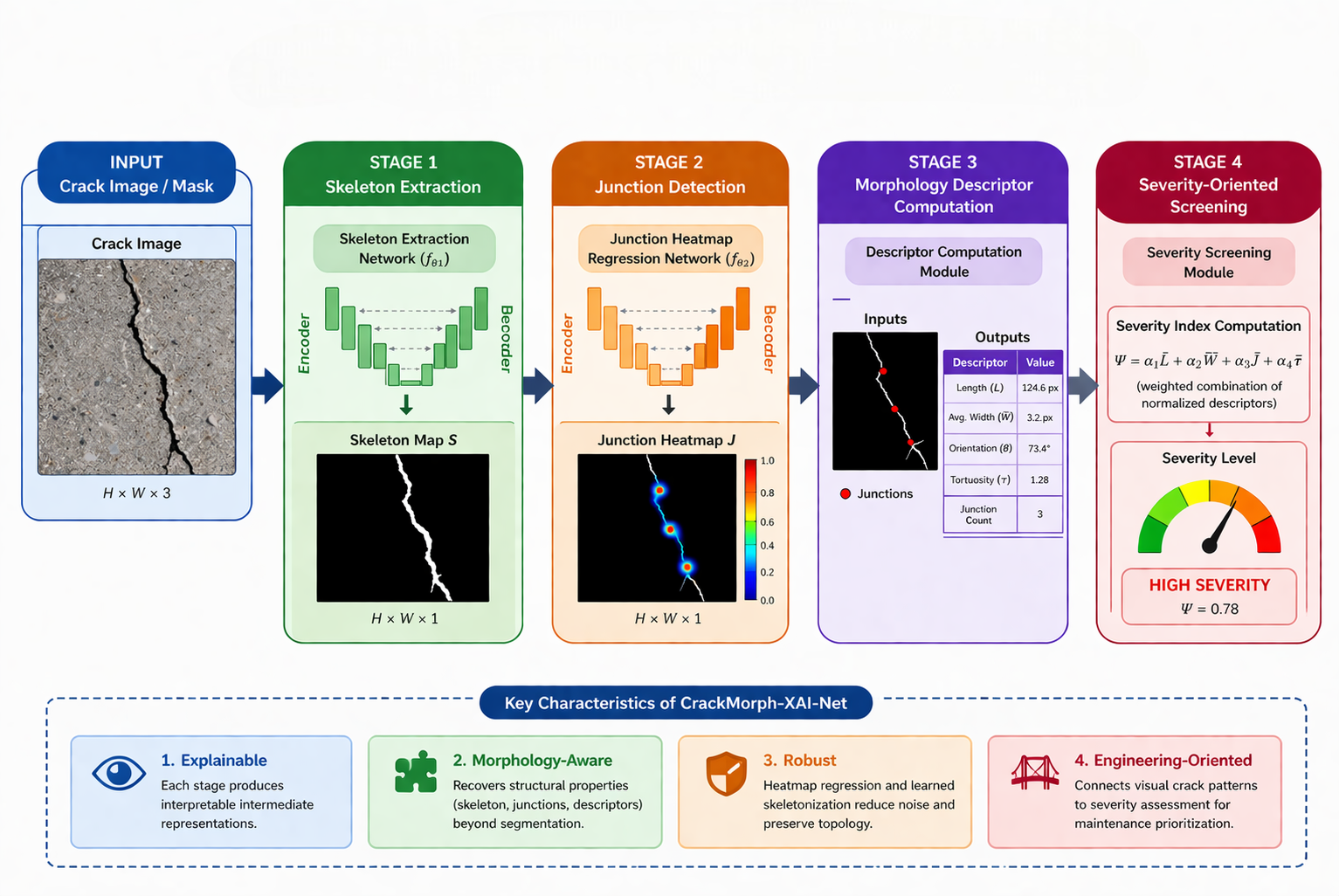}
\caption{Overall architectural structure of the proposed morphology-aware crack analysis framework. The sequential pipeline effectively translates raw crack-region data into continuous skeletons, probabilistic junctions, explicit morphology descriptors, and practical severity-oriented screening assessments.}
\label{fig:framework}
\end{figure}

Figure~\ref{fig:framework} provides a comprehensive schematic overview of the framework. Crucially, each stage is engineered to produce a distinct, interpretable intermediate output. The extracted skeleton delineates the continuous crack centerline, the junction heatmap probabilistically identifies regions of structural branching, the descriptor stage reports measurable geometric properties, and the final severity stage synthesizes these measurements into an actionable, screening-oriented interpretation.

A modular, cascade design is deliberately employed because the constituent subtasks fundamentally differ in their desired output modalities and network optimization behaviors. Skeleton extraction operates as a dense structural prediction task, whereas junction detection is inherently a sparse point-localization challenge. Furthermore, descriptor computation is a deterministic, geometric measurement operation, leading ultimately to a heuristic downstream interpretation task. By disentangling these processes, the cascade architecture allows each discrete stage to leverage task-appropriate supervision and objective functions while rigidly preserving transparency and verifiability throughout the computational pipeline.

\subsection{Skeleton Extraction}

The foundational stage of the pipeline addresses the recovery of a mathematically thin, one-pixel-wide topological skeleton from the broader crack region representation. An optimal skeletal representation must remain strictly centered within the physical crack trajectory, rigorously preserve longitudinal connectivity, resist artifactual fragmentation, and actively suppress the generation of spurious micro-branches. While classical morphological thinning algorithms can generate approximate skeletons, they are notoriously brittle and hypersensitive to noisy segmentation boundaries, often producing erratic, jagged centerlines. To achieve superior structural robustness, skeleton recovery is instead formulated as a fully learned structural prediction problem.

To execute this, a U-Net-style encoder--decoder network is deployed. The encoder pathway effectively captures multi-scale semantic context through successive downsampling and deep convolutional feature extraction, while the expansive decoder pathway reconstructs a spatially precise, full-resolution skeleton map. Symmetrical skip connections are utilized to concatenate high-resolution spatial feature maps from the encoder directly into the decoder, thereby mitigating spatial degradation and heavily aiding in the preservation of continuous, thread-like crack structures.

The composite loss function driving the skeleton extraction network seamlessly integrates Dice loss and binary cross-entropy (BCE):
\begin{equation}
\mathcal{L}_{skel} =
\mathcal{L}_{Dice}(S,\hat{S}) + \lambda_1 \mathcal{L}_{BCE}(S,\hat{S}),
\label{eq:skel_loss}
\end{equation}
where $S$ represents the ground-truth reference skeleton map, $\hat{S}$ denotes the predicted skeleton tensor, and $\lambda_1$ serves as a weighting hyperparameter controlling the relative contribution of the BCE term. Dice-based loss formulations are highly effective for optimizing against severely imbalanced foreground-background targets \citep{sudre2017dice}. In this context, the Dice term aggressively encourages global spatial overlap for the highly sparse centerline targets, while the BCE term simultaneously enforces stable, pixel-wise probabilistic learning.

The resulting predicted skeleton serves as a highly interpretable structural object. It facilitates immediate visual verification by end-users and provides the direct geometric foundation required for subsequent length, tortuosity, orientation, and complex topological measurements.

\subsection{Junction Detection}

Building upon the recovered centerline, the second stage focuses on the precise localization of crack junctions, representing critical branch points or structural intersections. The accurate identification of these junctions is paramount, as they directly dictate the macroscopic complexity and topological classification of the defect network. Formulating this as a direct pixel-wise binary classification task is inherently problematic; junction pixels occupy an infinitesimally small fraction of the overall image domain, exacerbating spatial imbalance far beyond typical dense prediction tasks \citep{lin2017focal}. To circumvent this severe imbalance and smooth the optimization landscape, the task is redefined as continuous heatmap regression. This strategy draws inspiration from advanced keypoint localization methodologies commonly utilized in human pose estimation, where discrete targets are relaxed into spatial probability distributions \citep{newell2016hourglass}.

Specifically, for each annotated discrete junction coordinate $(x_i,y_i)$, a two-dimensional Gaussian response is generated to represent the target presence of a junction. The cumulative target heatmap is mathematically defined as:
\begin{equation}
H(x,y) = \max_i \exp\left(
-\frac{(x-x_i)^2 + (y-y_i)^2}{2\sigma^2}
\right),
\label{eq:heatmap}
\end{equation}
where the variance parameter $\sigma$ explicitly controls the spatial spread, or radius of influence, of the Gaussian peak.

The dedicated junction detection network is trained to predict this continuous probability distribution, $\hat{H}$, optimized via a standard mean squared error objective:
\begin{equation}
\mathcal{L}_{junc} =
\frac{1}{HW}\sum_{x=1}^{H}\sum_{y=1}^{W}
\left(H(x,y)-\hat{H}(x,y)\right)^2.
\label{eq:junc_loss}
\end{equation}

During the inference phase, candidate junctions are extracted by identifying local spatial maxima within the predicted heatmap. A predefined confidence threshold, coupled with a non-maximum suppression algorithm, is subsequently applied to filter out weak, spurious, or redundantly clustered detections. The output heatmap provides profound explainability, offering a visual gradient of where the neural network interprets structural bifurcation to be most probable.

\subsection{Morphology Descriptor Computation}

Transitioning from predictive modeling to quantitative measurement, the third stage translates the generated structural representations into an array of deterministic morphological descriptors. Unlike the preceding neural network-driven stages, this operation is strictly geometric and analytical. It computes essential engineering metrics directly from the predicted skeleton, the boundary crack mask, and the localized junction coordinates.

The framework systematically computes the following structural descriptors:
\begin{enumerate}
    \item \textbf{Crack length:} Accurately estimated as the connectivity-aware, pixel-wise path length traced along the extracted skeleton.
    \item \textbf{Average width:} Approximated globally by evaluating the ratio between the total semantic crack area and the derived skeleton length.
    \item \textbf{Orientation:} Computed by analyzing the dominant principal axis of the localized spatial skeleton coordinates.
    \item \textbf{Junction count:} Quantified directly from the explicitly detected and filtered branch points.
    \item \textbf{Tortuosity:} Calculated as the ratio between the actual skeleton path length and the direct endpoint-to-endpoint distance.
    \item \textbf{Topology class:} A categorical descriptor derived directly from branching complexity and junction abundance.
\end{enumerate}

Mathematically, the average global crack width is approximated using the relation:
\begin{equation}
W_{avg} \approx \frac{A}{L},
\label{eq:width}
\end{equation}
where $A$ represents the total pixel area of the binary crack mask, and $L$ denotes the path length of the skeleton. Furthermore, tortuosity—a critical indicator of fracture unpredictability—is formally defined as:
\begin{equation}
T = \frac{L_{path}}{D_{end}},
\label{eq:tortuosity}
\end{equation}
where $L_{path}$ is the integrated skeleton path length, and $D_{end}$ represents the Euclidean, straight-line distance spanning the primary crack endpoints. Consequently, a higher tortuosity metric actively indicates a significantly more irregular, meandering, or mechanically deflected crack path.

Finally, global topology is categorized purely on junction-based structural complexity. Under this schema, \textit{linear} cracks possess no junctions, \textit{branched} cracks feature one or two localized junctions, \textit{complex} anomalies exhibit three to five junctions, and extensive \textit{network} cracks contain greater than five junctions. This structured categorization delivers a highly interpretable, macroscopic summary of the damage organization.

\begin{figure}[t]
\centering
\includegraphics[width=0.98\linewidth]{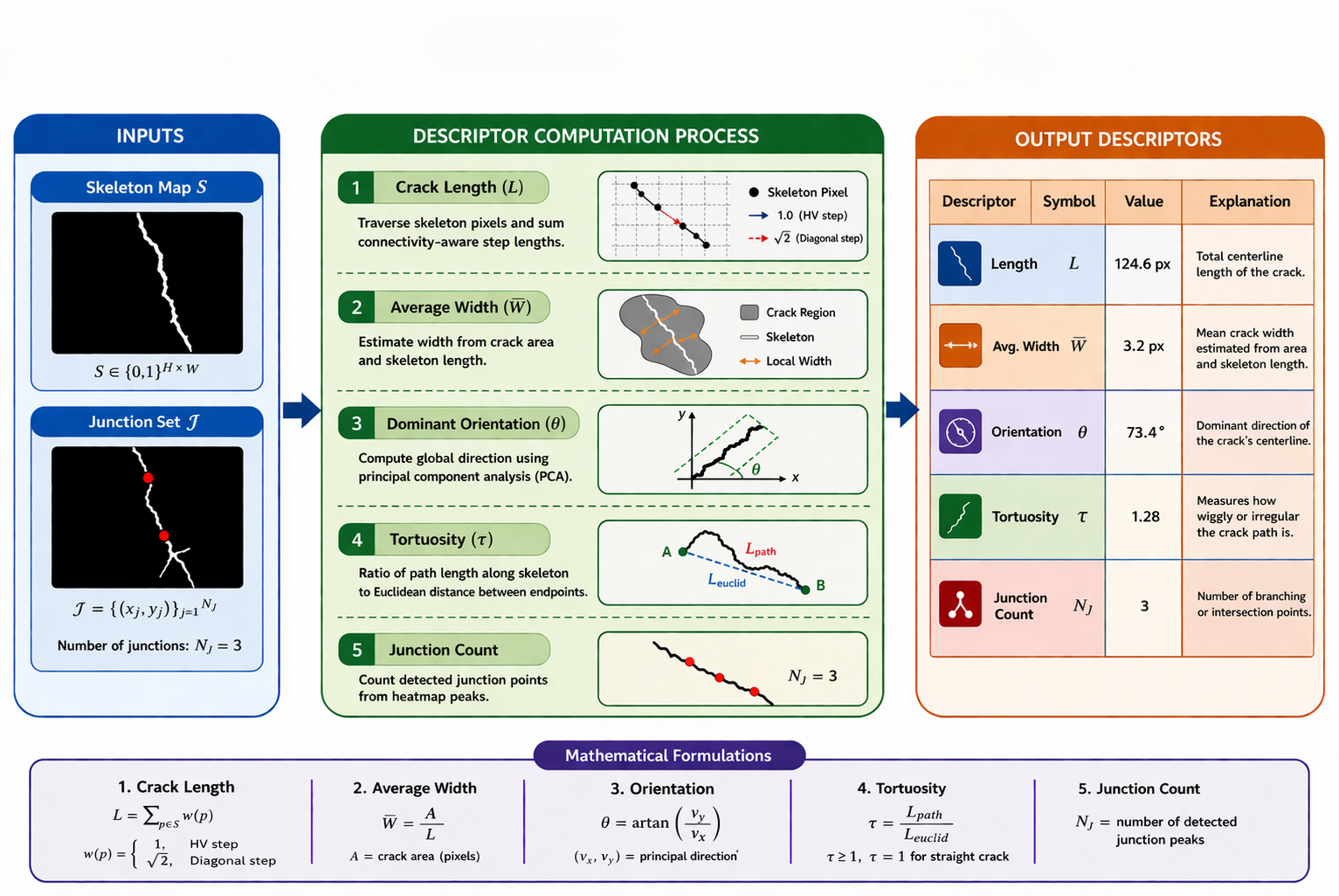}
\caption{Deterministic computation of morphology descriptors leveraging the predicted skeleton and mapped junction structure. This stage analytically bridges computer vision outputs with measurable, engineering-relevant geometric quantities.}
\label{fig:descriptors}
\end{figure}

\subsection{Severity-Oriented Screening}

The culmination of the pipeline yields a severity-oriented screening metric, synthesizing the computed morphological descriptors into an actionable assessment. It is crucial to note that this final stage is not intended to supplant rigorous, material-specific finite element analysis or in situ fracture mechanics evaluations. Rather, it functions as a highly efficient, data-driven prioritization heuristic aimed at triaging infrastructure damage and flagging critical defects for human expert review.

The foundational logic of this screening stage is inspired by the established principles of linear elastic fracture mechanics (LEFM), which state that absolute crack severity is fundamentally dictated by a combination of geometric bounds and structural complexity (e.g., length, width, and localized branching). This physical behavior is frequently abstracted through the Mode-I stress intensity factor \citep{williams1957stress,anderson2017fracture,mallikarjunaiah2025hp,Mallikarjunaiah2026,Gou2015}:
\begin{equation}
K_I = Y\sigma \sqrt{\pi a},
\label{eq:sif}
\end{equation}
where $\sigma$ denotes the far-field applied stress, $a$ serves as a governing parameter intrinsically related to crack length, and $Y$ functions as a dimensionless geometry factor that accounts for structural boundary conditions. Within the proposed framework, this classical relationship is utilized conceptually to inform the severity weighting, rather than attempting a full, computationally prohibitive simulation of material-specific failure. The resultant severity output should therefore be treated strictly as a screening-oriented, decision-support indicator, designed to accelerate maintenance planning rather than serving as a definitive declaration of structural safety.

\section{Experimental Setup}
\label{sec:experiments}

\subsection{Training Protocol}

To ensure optimal convergence and prevent the network optimization pathways from interfering with one another, the skeleton extraction and junction detection stages are trained independently utilizing the augmented CRACK500-Extended annotations. Prior to network ingestion, all images undergo a standardized preprocessing pipeline where they are uniformly resized to a spatial resolution of $640 \times 640$ pixels, converted to single-channel grayscale to minimize computational overhead, and intensity-normalized to the continuous range of $[0,1]$. 

The skeleton extraction network is optimized using the composite Dice and binary cross-entropy loss function defined in Eq.~\eqref{eq:skel_loss}, which simultaneously mitigates the extreme class imbalance of centerline targets while stabilizing pixel-wise probabilistic learning. Conversely, the junction detection network is optimized using the continuous heatmap regression loss outlined in Eq.~\eqref{eq:junc_loss}, smoothing the penalty landscape for highly sparse structural branch points. Both networks utilize the Adam optimizer \citep{kingma2015adam} with an initial learning rate fixed at $1\times10^{-4}$. To rigorously guard against overfitting, validation-based model selection is employed, dynamically halting training and preserving the model weights that achieve the lowest validation loss across a defined patience window. The complete experimental configuration is summarized in Table~\ref{tab:setup}.

\begin{table}[htpb]
\centering
\caption{Summary of the experimental configuration and hyperparameter settings.}
\label{tab:setup}
\begin{tabular}{ll}
\toprule
Item & Setting \\
\midrule
Base dataset & CRACK500 \\
Extended annotations & Skeleton maps, junction heatmaps, topology labels \\
Train/validation/test split & 250 / 50 / 200 images \\
Augmented training samples & 1,896 \\
Input resolution & $640 \times 640$ \\
Input format & Grayscale, normalized to $[0,1]$ \\
Skeleton model & U-Net-style encoder--decoder \\
Junction model & U-Net-style heatmap regressor \\
Skeleton loss & Dice + binary cross-entropy \\
Junction loss & Mean squared error \\
Optimizer & Adam \\
Initial learning rate & $1\times10^{-4}$ \\
\bottomrule
\end{tabular}
\end{table}

\subsection{Algorithmic Implementation}

To facilitate reproducibility and explicitly define the computational flow of CrackMorph-XAI-Net, the procedures for both the training phase and the subsequent inference and descriptor extraction phase are formalized below. Algorithm~\ref{alg:training} delineates the independent optimization of the structural subnetworks, while Algorithm~\ref{alg:inference} details the deterministic, end-to-end pipeline utilized during practical evaluation and field deployment.

\begin{algorithm}[htpb]
\caption{CrackMorph-XAI-Net: Training Protocol}
\label{alg:training}
\begin{algorithmic}[1]
\REQUIRE Augmented training set $\mathcal{D}_{train} = \{(I_i, S_i, H_i)\}_{i=1}^N$, Learning rate $\alpha$, Max epochs $E$
\ENSURE Optimized network weights $\Theta_{skel}^*$, $\Theta_{junc}^*$
\STATE Initialize weights $\Theta_{skel}, \Theta_{junc}$ randomly or via pre-training
\FOR{epoch $= 1$ \TO $E$}
    \FOR{each mini-batch $b \subset \mathcal{D}_{train}$}
        \STATE \textit{// Stage 1: Skeleton Network Optimization}
        \STATE Predict skeleton probability map: $\hat{S}_b \leftarrow \text{SkeletonNet}(I_b; \Theta_{skel})$
        \STATE Compute composite loss $\mathcal{L}_{skel}$ using $S_b$ and $\hat{S}_b$ via Eq.~\eqref{eq:skel_loss}
        \STATE Update weights: $\Theta_{skel} \leftarrow \Theta_{skel} - \alpha \nabla_{\Theta} \mathcal{L}_{skel}$
        
        \STATE \textit{// Stage 2: Junction Network Optimization}
        \STATE Predict continuous heatmap: $\hat{H}_b \leftarrow \text{JunctionNet}(I_b; \Theta_{junc})$
        \STATE Compute regression loss $\mathcal{L}_{junc}$ using $H_b$ and $\hat{H}_b$ via Eq.~\eqref{eq:junc_loss}
        \STATE Update weights: $\Theta_{junc} \leftarrow \Theta_{junc} - \alpha \nabla_{\Theta} \mathcal{L}_{junc}$
    \ENDFOR
    \STATE Evaluate $\mathcal{L}_{skel}$ and $\mathcal{L}_{junc}$ on validation set $\mathcal{D}_{val}$
    \STATE Save $\Theta_{skel}^*$ and $\Theta_{junc}^*$ if validation loss improves (Early Stopping)
\ENDFOR
\end{algorithmic}
\end{algorithm}

\begin{algorithm}[htpb]
\caption{CrackMorph-XAI-Net: Inference and Morphology Extraction}
\label{alg:inference}
\begin{algorithmic}[1]
\REQUIRE Test image $I$, Crack region mask $M$, Trained weights $\Theta_{skel}^*$, $\Theta_{junc}^*$
\ENSURE Morphology descriptors $\mathcal{M}$, Topology Class $\mathcal{T}$, Severity Indicator
\STATE $I_{norm} \leftarrow \text{Preprocess}(I)$
\STATE \textit{// Execute Learned Structural Predictions}
\STATE $\hat{S}_{prob} \leftarrow \text{SkeletonNet}(I_{norm}; \Theta_{skel}^*)$
\STATE $\hat{S}_{binary} \leftarrow \text{Apply spatial threshold to } \hat{S}_{prob}$
\STATE $\hat{H}_{prob} \leftarrow \text{JunctionNet}(I_{norm}; \Theta_{junc}^*)$
\STATE $\mathcal{J} \leftarrow \text{Extract spatial coordinates of local maxima from } \hat{H}_{prob}$ via NMS
\STATE \textit{// Compute Deterministic Morphology Descriptors}
\STATE $L \leftarrow \text{Calculate path length along } \hat{S}_{binary}$
\STATE $W_{avg} \leftarrow \text{Area}(M) / L$ \hfill \textit{(Eq.~\ref{eq:width})}
\STATE $D_{end} \leftarrow \text{Euclidean distance between endpoints of } \hat{S}_{binary}$
\STATE $T \leftarrow L / D_{end}$ \hfill \textit{(Eq.~\ref{eq:tortuosity})}
\STATE Count total localized junctions: $N_J = |\mathcal{J}|$
\STATE Compile descriptor array: $\mathcal{M} \leftarrow \{L, W_{avg}, T, N_J, \text{Orientation}\}$
\STATE \textit{// Derive Hierarchical Classifications}
\STATE $\mathcal{T} \leftarrow \text{Classify topology based on } N_J \text{ (Linear, Branched, Complex, Network)}$
\STATE Compute screening-oriented Severity Indicator using $\mathcal{M}$ \hfill \textit{(Conceptual Eq.~\ref{eq:sif})}
\RETURN $\mathcal{M}, \mathcal{T}, \text{Severity Indicator}$
\end{algorithmic}
\end{algorithm}

\subsection{Evaluation Metrics}

To comprehensively assess the framework across its distinct operational stages, a multi-tiered evaluation strategy is employed. The skeleton extraction stage, being a dense structural prediction task, is quantitatively evaluated using the Dice coefficient, Intersection-over-Union (IoU), pixel-wise precision, and recall. Furthermore, topology preservation is uniquely measured to ensure the predicted centerline maintains structural fidelity without artefactual fragmentation. 

The junction detection stage is evaluated as a point-localization challenge using precision, recall, F1-score, average precision, and matched true/false positives. To account for minor spatial variances inherent in heatmap regression, a predicted junction is classified as a correct detection (true positive) if its coordinate lies within a strictly defined, acceptable spatial neighborhood radius of a ground-truth reference junction.

Moving to the deterministic measurement stage, descriptor accuracy is subjected to rigorous statistical evaluation using Mean Absolute Error (MAE), Root Mean Squared Error (RMSE), Mean Absolute Percentage Error (MAPE), and the Pearson correlation coefficient to demonstrate agreement with reference geometric values. Topology classification is assessed via overall accuracy and a detailed confusion matrix to illustrate class-wise prediction behaviors. 

Finally, to conclusively demonstrate the efficacy of the learned framework, benchmark comparisons are conducted against a classical, unlearned structural pipeline. This baseline consists of traditional morphological mask thinning, neighbor-count based junction extraction, and identical deterministic descriptor computation, thereby isolating and highlighting the performance gains achieved specifically through the proposed deep learning architectures.

\section{Results and Discussion}
\label{sec:results_discussion}

The fundamental premise of this work is that automated crack analysis must evolve beyond simple pixel-level localization to provide actionable, structural intelligence. While binary masks successfully delineate damage regions, they lack the intrinsic geometric organization required for rigorous engineering interpretation. The proposed CrackMorph-XAI-Net addresses this critical gap. This section details the quantitative evaluation of each progressive stage of the framework, immediately contextualizing these metrics with insights into architectural behavior, physical interpretability, and practical limitations.

\subsection{Learned Skeleton Extraction}

The foundational stage of the proposed pipeline is the recovery of a continuous, one-pixel-wide centerline representation. The integrity of this skeleton is paramount; if the centerline is fragmented or populated with artefactual micro-branches, all subsequent morphological computations become inherently unreliable. 

Table~\ref{tab:skeleton} details the performance of the learned skeleton extraction module on the CRACK500-Extended test set. The architecture achieves an exceptional mean Dice coefficient of 0.991 and a mean Intersection-over-Union (IoU) of 0.983. The tightly coupled precision and recall scores indicate that the model aggressively recovers ground-truth centerline pixels while successfully suppressing spurious structural noise. Crucially, macroscopic topology is strictly preserved in 98.5\% of the evaluated images (197 out of 200), demonstrating remarkable geometric fidelity.

\begin{table}[htpb]
\centering
\caption{Quantitative evaluation of learned skeleton extraction performance on the CRACK500-Extended test set.}
\label{tab:skeleton}
\begin{tabular}{lcccc}
\toprule
Metric & Mean & Std. Dev. & Median & Range \\
\midrule
Dice coefficient & 0.991 & 0.014 & 0.995 & [0.932, 1.000] \\
IoU & 0.983 & 0.025 & 0.990 & [0.872, 1.000] \\
Precision & 0.989 & 0.019 & 0.994 & [0.891, 1.000] \\
Recall & 0.994 & 0.011 & 0.997 & [0.945, 1.000] \\
Topology preservation & \multicolumn{4}{c}{98.5\% (197/200 images)} \\
\bottomrule
\end{tabular}
\end{table}

The necessity of formulating skeletonization as a learned structural prediction task rather than a post-processing heuristic is starkly illuminated when compared to classical methods. As evidenced in Table~\ref{tab:skeleton_baseline}, the proposed predictive approach dramatically eclipses the classical morphological thinning baseline, elevating the skeleton Dice score from 0.912 to 0.991. More importantly, topology preservation leaps by 11.5 percentage points. Traditional thinning algorithms remain notoriously hypersensitive to irregular mask boundaries, frequently generating jagged, physically implausible centerlines. By directly learning the structural mapping from aligned targets, the proposed framework effectively inoculates the skeleton against boundary-induced artifacts.

\begin{figure}[htpb]
\centering
\includegraphics[width=0.98\linewidth]{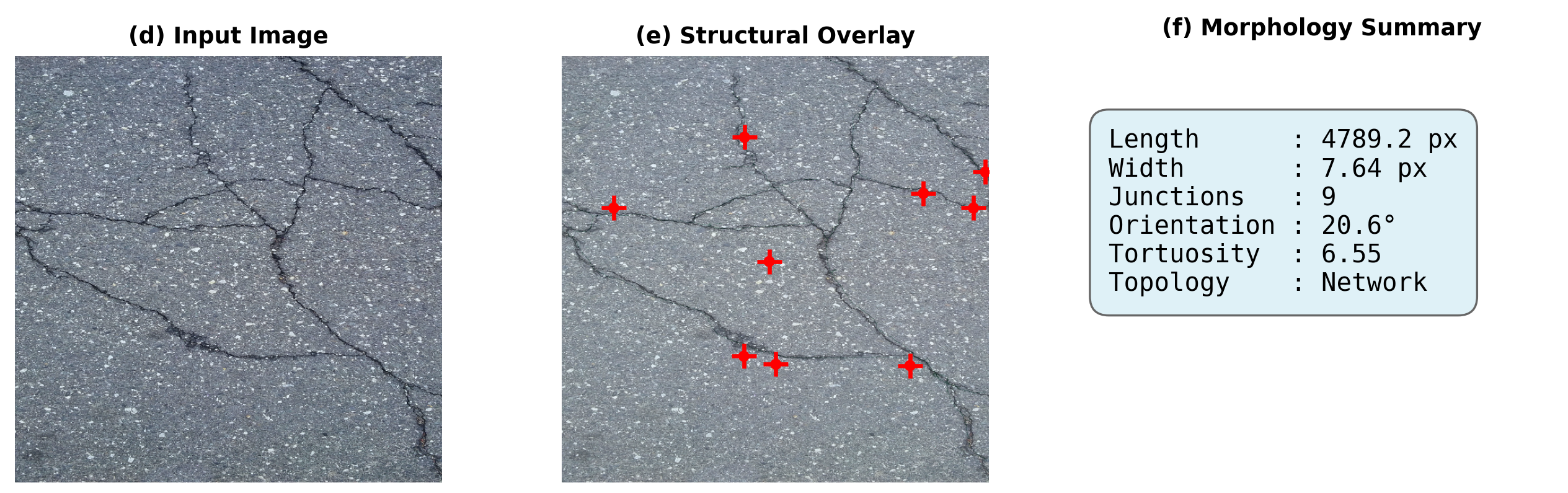}
\caption{Representative structural overlay and morphology summary for a test instance. The reliably recovered centerline and localized branch coordinates provide an interpretable foundation for subsequent descriptor computation.}
\label{fig:skeleton_results}
\end{figure}

\begin{table}[htpb]
\centering
\caption{Performance comparison between the proposed learned skeleton extraction and the classical thinning baseline.}
\label{tab:skeleton_baseline}
\begin{tabular}{lcc}
\toprule
Method & Skeleton Dice & Topology preservation \\
\midrule
Classical thinning baseline & $0.912 \pm 0.067$ & 87.0\% \\
Proposed skeleton extraction & $0.991 \pm 0.014$ & 98.5\% \\
Improvement & +0.079 & +11.5 percentage points \\
\bottomrule
\end{tabular}
\end{table}

\subsection{Continuous Junction Detection}

Identifying structural intersections—or junctions—is critical for defining the macroscopic complexity of a crack network. However, the inherent sparsity of these junctions presents a formidable optimization hurdle; binary classification of such minuscule pixel targets typically results in severe class imbalance and poor model convergence. 

To circumvent this, the framework deploys a Gaussian heatmap regression strategy. As reported in Table~\ref{tab:junction}, this approach yields a recall of 0.964 and an F1-score of 0.887. Across the entire test corpus, the model successfully localizes 301 out of 312 reference junctions, missing a mere 11 instances. While the precision score of 0.822 indicates the occasional generation of redundant candidate junctions, maintaining an aggressively high recall is strategically preferred for screening-oriented structural analysis. Missing a genuine structural branch exerts a disproportionately negative distortion on the overarching topology classification, whereas false positives are more easily filtered during deterministic descriptor evaluation.

\begin{table}[htpb]
\centering
\caption{Junction detection performance utilizing Gaussian heatmap regression on the test set.}
\label{tab:junction}
\begin{tabular}{lccc}
\toprule
Metric & Value & 95\% CI & Per-image mean \\
\midrule
Precision & 0.822 & [0.784, 0.860] & 0.831 \\
Recall & 0.964 & [0.941, 0.987] & 0.971 \\
F1-score & 0.887 & [0.862, 0.912] & 0.893 \\
Average precision & 0.746 & [0.701, 0.791] & -- \\
True positives & 301 & -- & 1.51 \\
False positives & 65 & -- & 0.33 \\
False negatives & 11 & -- & 0.06 \\
\bottomrule
\end{tabular}
\end{table}

The specific advantage of heatmap regression over direct binary prediction is quantified in Table~\ref{tab:heatmap}. By relaxing the discrete point target into a continuous spatial probability distribution, the regression model improves recall by 0.093 and the overall F1-score by 0.057. This smoother optimization landscape forces the network to learn robust local structural context rather than merely memorizing isolated pixel intensities.

\begin{figure}[htpb]
\centering
\includegraphics[width=0.98\linewidth]{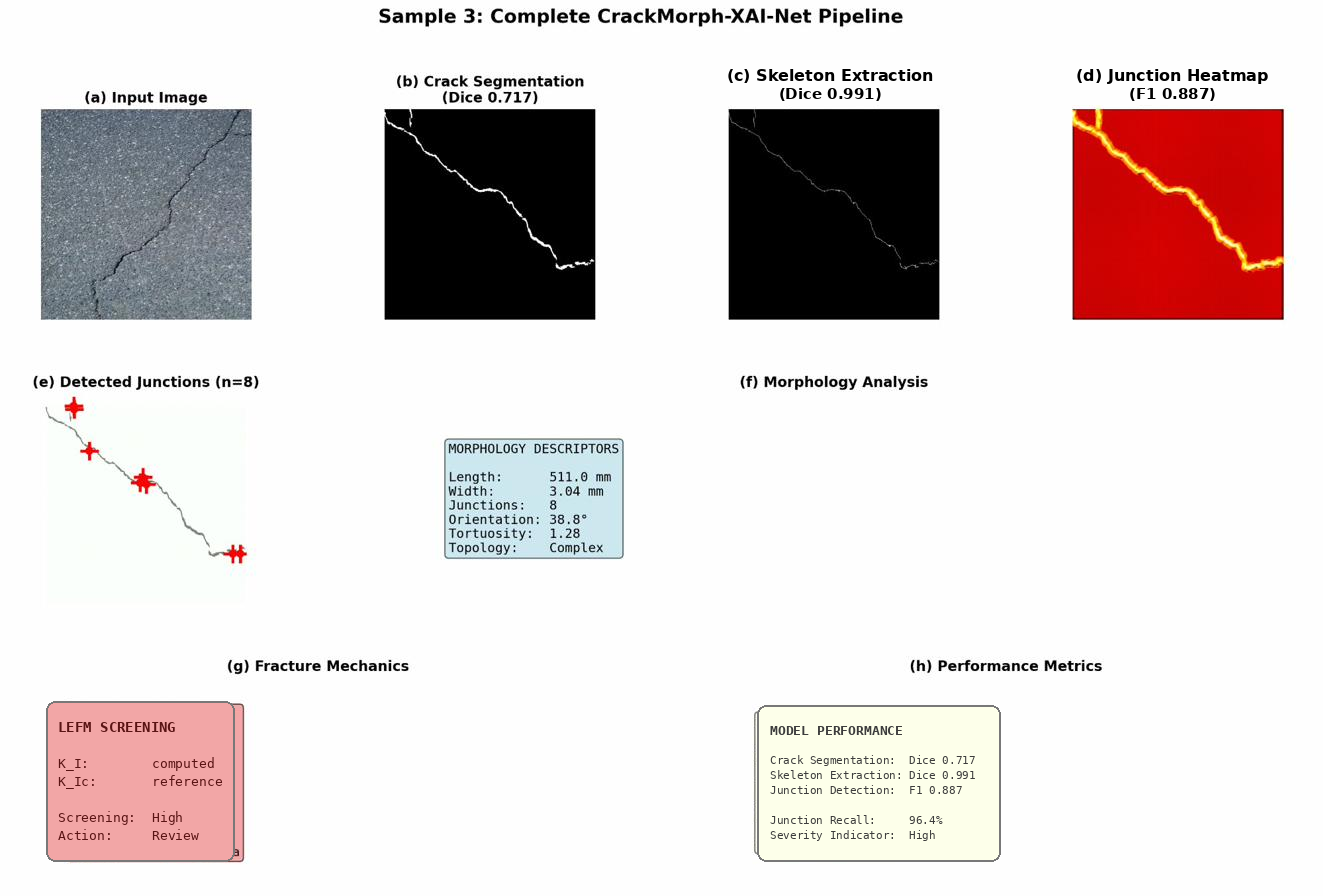}
\caption{A comprehensive visualization of the sequential pipeline, transitioning from raw image segmentation through learned skeletonization and continuous heatmap prediction, culminating in fully interpretable morphology measurements.}
\label{fig:junction_results}
\end{figure}

\begin{table}[htpb]
\centering
\caption{Ablative comparison demonstrating the superiority of Gaussian heatmap regression over standard binary junction classification.}
\label{tab:heatmap}
\begin{tabular}{lccc}
\toprule
Method & Precision & Recall & F1-score \\
\midrule
Binary junction classification & 0.793 & 0.871 & 0.830 \\
Gaussian heatmap regression & 0.822 & 0.964 & 0.887 \\
Improvement & +0.029 & +0.093 & +0.057 \\
\bottomrule
\end{tabular}
\end{table}

\subsection{Morphology Measurement and Topology Classification}

The ultimate utility of the CrackMorph-XAI-Net framework rests on its capacity to translate intermediate visual outputs into rigorous, engineering-grade numerical measurements. Table~\ref{tab:descriptors} illustrates profound agreement between the predicted morphological descriptors and ground-truth reference values. Notably, the Pearson correlation coefficients consistently exceed 0.95 across all metrics, with physical path length and axis orientation peaking at 0.982 and 0.991, respectively. This confirms that the learned structural representations are not merely visually plausible, but mathematically reliable enough to support quantitative health monitoring protocols. Because these metrics are derived deterministically from visible intermediate states, their derivations remain fully traceable—cementing the framework's claim to intrinsic explainability.

\begin{table}[htpb]
\centering
\caption{Accuracy of the deterministically computed morphology descriptors, illustrating strong agreement with reference measurements.}
\label{tab:descriptors}
\begin{tabular}{lcccc}
\toprule
Descriptor & MAE & RMSE & MAPE (\%) & Pearson $r$ \\
\midrule
Length (mm) & 12.3 & 18.7 & 8.4 & 0.982 \\
Width (mm) & 0.31 & 0.49 & 11.2 & 0.957 \\
Orientation ($^\circ$) & 3.8 & 6.2 & -- & 0.991 \\
Junction count & 0.18 & 0.44 & 9.5 & 0.988 \\
Tortuosity & 0.042 & 0.067 & 3.6 & 0.974 \\
\bottomrule
\end{tabular}
\end{table}

In parallel, topology classification abstracts the localized junction counts into a macroscopic summary of branching severity. As shown in the confusion matrix in Table~\ref{tab:topology}, the system achieves an overall categorical accuracy of 83.5\%. An analysis of the misclassifications reveals that errors are almost exclusively localized between adjacent complexity boundaries (e.g., distinguishing a heavily branched crack from a complex network). This behavior suggests that even when exact branch counting encounters ambiguity in highly degraded regions, the framework reliably preserves the general magnitude of structural distress.

\begin{table}[htpb]
\centering
\caption{Confusion matrix detailing topology classification performance across four distinct complexity thresholds.}
\label{tab:topology}
\begin{tabular}{lccccc}
\toprule
True / Predicted & Linear & Branched & Complex & Network & Total \\
\midrule
Linear & 110 & 4 & 0 & 0 & 114 \\
Branched & 8 & 45 & 4 & 0 & 57 \\
Complex & 0 & 6 & 16 & 1 & 23 \\
Network & 0 & 0 & 1 & 5 & 6 \\
\midrule
Overall accuracy & \multicolumn{5}{c}{83.5\%} \\
\bottomrule
\end{tabular}
\end{table}

\subsection{Comprehensive Baseline Analysis and Architecture Ablation}

To rigorously validate the overall architectural design, Table~\ref{tab:baseline} compares the complete CrackMorph-XAI-Net pipeline against a fully classical structural baseline (comprising masking, algorithmic thinning, and neighbor-count junction extraction). The proposed deep learning framework dominates across every evaluated metric. The most transformative improvements are observed in the detection of junctions—where recall surges by 0.181—and in the drastic reduction of downstream measurement errors, dropping length MAE by 37.6\% and width MAE by 35.4\%. These combined metrics unequivocally establish that reliance on classical post-processing of segmentation masks is inherently brittle, whereas stage-wise structural learning yields highly resilient engineering profiles.

\begin{table}[htpb]
\centering
\caption{Comprehensive end-to-end comparison between the classical structural analysis pipeline and the proposed learned framework.}
\label{tab:baseline}
\begin{tabular}{lccc}
\toprule
Metric & Baseline & Proposed framework & Improvement \\
\midrule
Skeleton Dice & $0.912 \pm 0.067$ & $0.991 \pm 0.014$ & +0.079 \\
Topology preservation & 87.0\% & 98.5\% & +11.5 pp \\
Junction recall & 0.783 & 0.964 & +0.181 \\
Junction F1-score & 0.711 & 0.887 & +0.176 \\
Length MAE (mm) & 19.7 & 12.3 & -37.6\% \\
Width MAE (mm) & 0.48 & 0.31 & -35.4\% \\
Topology accuracy & 76.5\% & 83.5\% & +7.0 pp \\
\bottomrule
\end{tabular}
\end{table}

Furthermore, the architectural ablation study detailed in Table~\ref{tab:ablation} justifies the deployment of an independent cascade framework rather than a monolithic multi-task network. While multi-task architectures are computationally economical, forcing the network to simultaneously optimize dense centerline recovery alongside highly sparse junction localization creates detrimental gradient interference. The isolated cascade approach proves superior, significantly bolstering junction recall and overall topology accuracy. Additionally, the study confirms that aggressive spatial augmentation is necessary to simulate the immense structural diversity required to properly generalize these models.

\begin{table}[htpb]
\centering
\caption{Ablation study isolating the impact of major architectural and experimental design choices.}
\label{tab:ablation}
\begin{tabular}{lccc}
\toprule
Configuration & Skeleton Dice & Junction measure & Topology accuracy \\
\midrule
Unified multi-task model & $0.987 \pm 0.019$ & $0.894 \pm 0.143$ (recall) & 78.5\% \\
Cascade framework & $0.991 \pm 0.014$ & $0.964 \pm 0.096$ (recall) & 83.5\% \\
\midrule
No augmentation & $0.981 \pm 0.023$ & $0.861 \pm 0.138$ (F1) & -- \\
With augmentation & $0.991 \pm 0.014$ & $0.887 \pm 0.119$ (F1) & -- \\
\bottomrule
\end{tabular}
\end{table}

\subsection{Interpretability, Practical Applications, and Limitations}

\begin{figure}[htpb]
\centering
\includegraphics[width=0.98\linewidth]{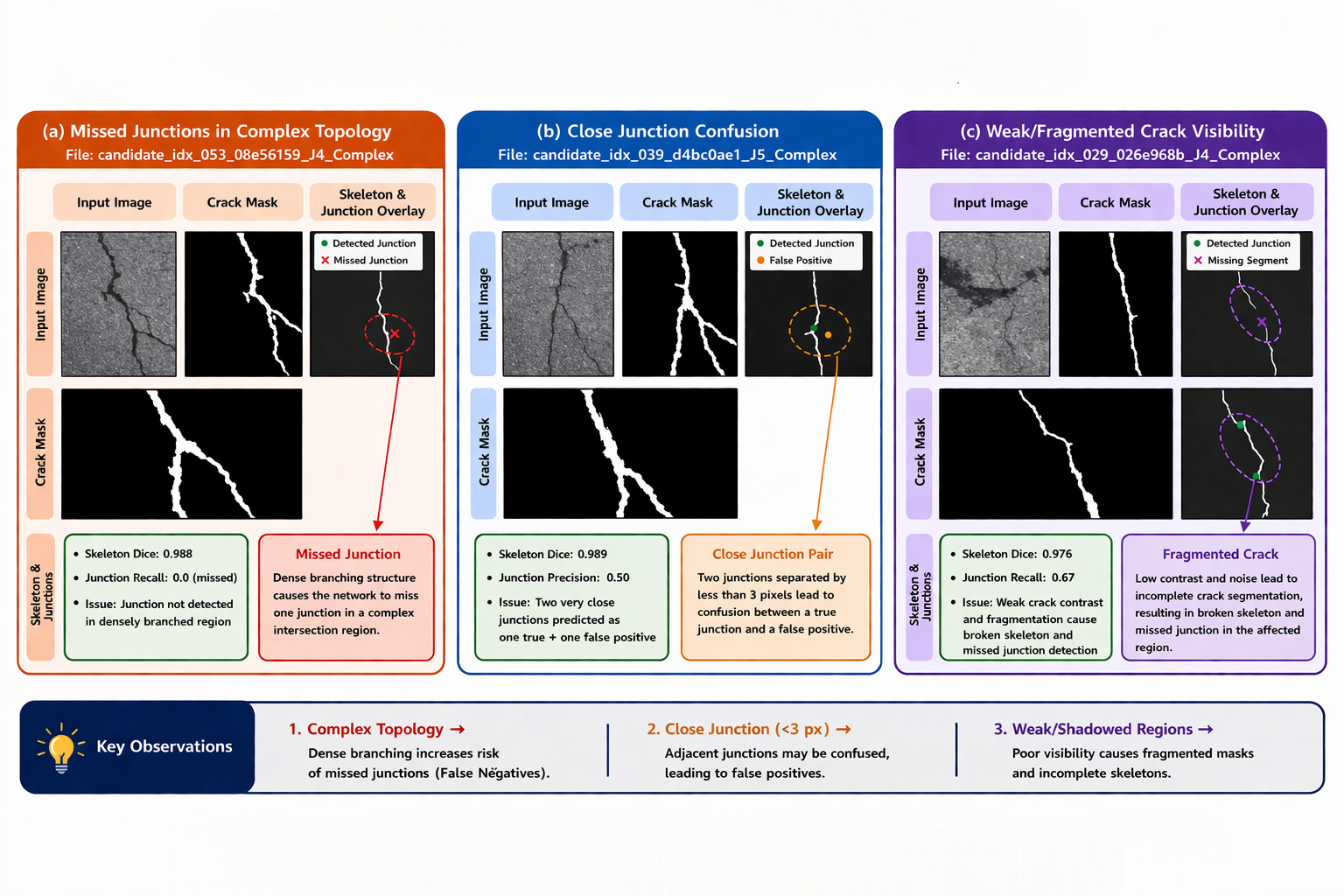}
\caption{Qualitative failure case analysis. The transparent nature of the pipeline allows engineers to explicitly diagnose prediction errors, such as merged junctions in highly complex topologies or fragmented centerlines in weak-contrast regions.}
\label{fig:qualitative}
\end{figure}

Beyond quantitative supremacy, the intrinsic value of CrackMorph-XAI-Net lies in its transparent explainability. Figure~\ref{fig:qualitative} provides qualitative examples of the framework's failure modes. Unlike black-box segmentation models where an erroneous prediction is an opaque dead end, the cascade design localizes the exact point of failure. If the final length calculation is inaccurate, an engineer can explicitly review the skeleton overlay to observe whether the centerline mistakenly traversed a region of dense noise. Errors typically manifest in low-contrast zones, exceptionally dense intersections, or closely spaced bifurcation points where Gaussian heatmaps may artificially merge. Because these intermediate states are visually accessible, domain experts can rapidly diagnose algorithmic shortcomings and manually intervene if critical structural details are obscured.

The ultimate severity-oriented screening output must be interpreted with clear engineering boundaries. While derived from the geometric parameters foundational to linear elastic fracture mechanics, this output is explicitly designed as a rapid triage heuristic, not a substitute for full-scale material analysis. True structural jeopardy is heavily governed by material heterogeneity, multi-axial loading states, and 3D boundary conditions—factors inherently absent from 2D visual inspections. Thus, the tool should be deployed as a decision-support mechanism to intelligently prioritize anomalous defects for human expert review.

Finally, several avenues for future refinement exist. The current empirical validation is constrained to macroscopic pavement imagery; rigorously validating the model's generalizability across diverse media—such as steel infrastructure or localized microstructural concrete spalling—remains an active necessity. Furthermore, the framework currently analyzes planar 2D morphology. Integrating multi-view stereopsis or laser profilometry to dynamically capture depth metrics would significantly elevate the fidelity of the severity screening protocol. Despite these limitations, this methodology marks a definitive, necessary paradigm shift for automated infrastructure inspection, demonstrating that algorithmic output must not only localize damage, but quantify its structural architecture in an inherently verifiable manner.

\section{Conclusion}
\label{sec:conclusion}

This paper introduced CrackMorph-XAI-Net, an end-to-end, intrinsically explainable framework designed to fundamentally advance the domain of automated infrastructure inspection. Moving beyond the historically restrictive paradigm of binary, pixel-level segmentation, the proposed architecture systematically maps raw visual data into a transparent hierarchy of structural representations. By sequentially generating topology-preserving skeletons, probabilistically mapped junction heatmaps, deterministic morphology descriptors, and severity-oriented heuristics, the framework successfully bridges the critical gap between computer vision outputs and actionable engineering intelligence. To rigorously validate this stage-wise methodology and catalyze future morphology-aware research, this study also developed CRACK500-Extended, a comprehensively enriched benchmark featuring strictly aligned spatial and topological annotations.

Extensive empirical evaluations unequivocally demonstrate the framework's quantitative and qualitative superiority over classical algorithmic post-processing. By formulating skeletonization as a learned structural prediction task, the model achieved a remarkable mean Dice coefficient of 0.991 and maintained absolute topological fidelity in 98.5\% of test instances, successfully eradicating the boundary-induced artifacts that notoriously plague traditional morphological thinning methods. Furthermore, resolving the extreme spatial sparsity of structural intersections via continuous Gaussian heatmap regression yielded an exceptional junction recall of 0.964 alongside an F1-score of 0.887. This robust sub-structural localization directly translates into highly reliable downstream geometric measurements; Pearson correlations exceeded 0.95 across all morphological descriptors, confirming that the framework's learned visual representations are mathematically consistent with true physical morphology.

Ultimately, the most significant contribution of CrackMorph-XAI-Net lies in its foundational commitment to structural transparency. By exposing every intermediate predictive state---from continuous centerlines to localized branch coordinates---the pipeline guarantees that final automated assessments are fully traceable. This provides domain experts with the explicit visual evidence necessary to trust, verify, or surgically diagnose the model's geometric reasoning, a critical requirement for deploying AI in safety-critical environments. Moving forward, research will prioritize evaluating the framework's robust generalizability across diverse material domains, such as steel infrastructure and complex concrete microstructures. Additionally, future iterations will explore the integration of spatiotemporal graph neural networks for tracking dynamic crack propagation, alongside tighter coupling with physics-informed material models to realize fully automated, predictive structural health monitoring.

\section*{Data and Code Availability}

The dataset extension and implementation details can be made available upon reasonable request, subject to the licensing and redistribution terms of the original CRACK500 dataset. 

\section*{Conflict of Interest}

The author declares no known competing financial interests or personal relationships that could have influenced the work reported in this paper.

%

\bibliographystyle{unsrtnat}

\end{document}